\documentclass[10pt]{amsart}
\usepackage{amssymb,enumerate}
\usepackage{amsmath}
\usepackage{graphics}
 \newcommand{\eR}{\mathbb R}
\newcommand{\E}{\mathbb E}
\newcommand{\rf}[1]{{\rm(\ref{#1})}}

\newtheorem{theorem}{Theorem}
\newtheorem{definition}{Definition}
\newtheorem{lemma}{Lemma}
\newtheorem{example}{Example}
\newtheorem{remark}{Remark}

\begin{document}

\title{Levin Ste\v{c}kin theorem and inequalities of the Hermite-Hadamard type} 
\author{Tomasz Szostok}
\address{Institute of Mathematics, University of Silesia, Bankowa 14, 40-007 Katowice, Poland}
\email{tszostok@math.us.edu.pl}

\begin{abstract} Recently Ohlin lemma on convex stochastic ordering was used to obtain some inequalities
of Hermite-Hadamard type. Continuing this idea, we use Levin-Ste\v{c}kin result to determine all inequalities of the forms:
$$\sum_{i=1}^3a_if(\alpha_ix+(1-\alpha_i)y)\leq \frac{1}{y-x}\int_{x}^yf(t),$$
$$a_1f(x)+\sum_{i=2}^3a_if(\alpha_ix+(1-\alpha_i)y)+a_4f(y)\geq \frac{1}{y-x}\int_{x}^yf(t)$$
and
$$af(\alpha_1x+(1-\alpha_1)y)+(1-a)f(\alpha_2x+(1-\alpha_2)y)\leq b_1f(x)+b_2f(\beta x+(1-\beta)y)+b_3f(y)$$
which are satisfied by all convex functions $f:[x,y]\to\eR.$
As it is easy to see, the same methods may be applied to deal
with longer expressions of the forms considered. As particular cases of our results we obtain 
some known inequalities.
\end{abstract}
\maketitle
\pagestyle{myheadings}
\section{Introduction}

In this paper we obtain some class of  inequalities of the Hermite-Hadamard type. 
First we write the classical Hermite-Hadamard inequality    
\begin{equation}
\label{HH}
f\left(\frac{x+y}{2}\right)\leq\frac{1}{y-x}\int_{x}^yf(t)dt\leq\frac{f(x)+f(y)}{2}
\end{equation}
(see \cite{DP} for many generalizations and applications of \rf{HH}).

In recent papers \cite{Rajba} and \cite{Szostok}  Ohlin lemma on convex stochastic ordering was used to obtain inequalities of the Hermite-Hadamard type. Therefore we cite this lemma.
\begin{lemma} (Ohlin \cite{Ohlin})
Let $X_1,X_2$ be two random variables such that $\E X_1=\E X_2$  and let  $F_1,F_2$ be their distribution functions.
If $F_1,F_2$ satisfy for some $x_0$ the following inequalities
\begin{equation}
\label{x0}
F_1(x)\leq F_2(x)\;{\rm if}\;x<x_0\;\;{\rm and}\;\;F_1(x)\geq F_2(x)\;{\rm if}\;x>x_0 
\end{equation}
then 
\begin{equation}
\label{m}
\E f(X_1)\leq \E f(X_2)
\end{equation}
 for all continuous and convex functions $f:\eR\to\eR.$
 \label{OL}
\end{lemma}

In paper \cite{Szostok} all numbers $a,\alpha,\beta\in[0,1]$ such that for all convex functions 
$f$ the inequality 
$$af(\alpha x+(1-\alpha )y)+(1-a)f(\beta x+(1-\beta) y)\leq \frac{1}{y-x}\int_{x}^yf(t)dt$$
is satisfied and all $a,b,c,\alpha\in(0,1)$ with $a+b+c=1$ for which we have
$$bf(x)+cf(\alpha x+(1-\alpha)y)+df(y)\geq\frac{1}{y-x}\int_{x}^y f(t)dt.$$ This means that 
the expressions   $f\left(\frac{x+y}{2}\right)$ and $\frac{f(x)+f(y)}{2}$ were replaced by more complicated ones.
As  particular cases, it was possible to obtain some known inequalities. We shall discuss these 
applications later in the paper.

However, the method used in \cite{Szostok} was valid only for these specific situations. Using
results proved there, it was not possible to get inequalities of the same  type for longer sums. Moreover, 
it was not even possible to obtain inequalities of the form
$$af(\alpha x+(1-\alpha )y)+(1-a)f(\beta x+(1-\beta y))\leq
bf(x)+cf(\alpha x+(1-\alpha)y)+df(y).$$
Therefore, in order to extend results from \cite{Szostok},  in the present approach we are going to use a result from \cite{LS}, (see also \cite{NP} Theorem 4.2.7).

\begin{theorem}(Levin, Ste\v{c}kin)
\label{LS}
Let $x,y\in\eR,x<y$ and let $F_1,F_2:[x,y]\to\eR$ be two functions with bounded variation such that $F_1(x)=F_2(x).$
Then, in order that 
$$\int_x^y f(t)dF_1(t)\leq \int_x^y f(t)dF_2(t),$$
it is necessary and sufficient that $F_1$ and $F_2$ verify the following three conditions:
\begin{equation}F_1(y)=F_2(y),\end{equation}
\begin{equation}\int_x^s F_1(t)dt\leq\int_x^s F_2(t)dt,\;s\in(x,y)
\label{FGi}
\end{equation}
and
\begin{equation}\int_x^y F_1(t)dt=\int_x^y F_2(t)dt.
\label{FGe}
\end{equation}
\end{theorem}

\begin{remark}Observe that if measures $\mu_1,\mu_2$ corresponding to the random variables
occurring in Ohlin lemma are concentrated on the interval $[x,y]$ then Ohlin lemma  is an easy consequence of Theorem \ref{LS}. Indeed, $\mu_1,\mu_2$ are probabilistic measures thus we have
$F_1(x)=F_2(x)=0$ and $F_1(y)=F_2(y)=1.$ Moreover $\E X_1=\E X_2$ yields \rf{FGe} and from the inequalities \rf{x0} we get \rf{FGi}.  
 \end{remark}

Now we shall use Theorem \ref{LS} to make an observation which is more general than 
Ohlin lemma and concerns the situation when functions $F_1,F_2$ have more crossing points than one.
First we need the following definition.
\begin{definition}
Let $F_1,F_2:[x,y]\to\eR$ be functions and let $x=x_0<x_1<\dots<x_n<x_{n+1}=y.$ We say that the pair 
$(F_1,F_2)$ crosses $n$-times (at points $x_1,\dots,x_n$) if the inequalities 
\begin{equation}
\label{xi}
F_1(t)\leq F_2(t),\;t\in(x_i,x_{i+1})\;{\rm and}\;F_1(t)\geq F_2(t),\;t\in(x_{i+1},x_{i+2})
\end{equation}
where $i$ is even, are satisfied and
$$\int_{x_i}^{x_{i+1}}F_1(t)-F_2(t)dt\neq 0$$
for all $i\in\{0,\dots,n\}.$
\end{definition}

In the next part of the paper we shall use a lemma which may be found in \cite{OS1}.
It is easy to prove this lemma using Theorem \ref{LS}. 
 For the sake of simplicity we omit the proof of this result.

\begin{lemma}
\label{mainl}
Let $F_1,F_2:[x,y]\to\eR$ be two functions with bounded variation such that $F_1(x)=F_2(x)$ and $F_1(y)=F_2(y)$
let $x_1,\dots,x_n\in(x,y)$ and let $(F_1,F_2)$ cross $n-$times, at $x_1,\dots,x_n\in(x,y).$

(i) If $n$ is even then the inequality
\begin{equation}
\label{maini}
\int_x^y f(t)dF_1(t)\leq \int_x^y f(t)dF_2(t)
\end{equation}
  is not satisfied by all convex functions $f:[x,y]\to\eR.$
  
(ii) If $n$ is odd then we define numbers $A_i$ by the following formula
$$A_i:=\int_{x_i}^{x_{i+1}}|F_1(x)-F_2(x)|dx.$$ 

Inequality \rf{maini} is satisfied for all convex functions $f:[x,y]\to\eR$  if and only if the following
inequalities hold true
$$A_0\geq A_1,$$
$$A_0-A_1+A_2\geq A_3,$$
$$A_0-A_1+A_2-A_3+A_4\geq A_5,$$
$$\vdots$$
$$A_0-A_1+A_2-A_3+A_4-A_5+\dots-A_{n-4}+A_{n-3}\geq A_{n-2}.$$
\end{lemma} 

\section{Results}

In this part of the paper we shall show how Lemma \ref{mainl} may be used to obtain 
inequalities of the Hermite-Hadamard type and inequalities between quadrature operators. As it is easy to 
see, if an inequality of this kind is satisfied for every convex function defined on the interval $[0,1]$
then it is satisfied by every convex function defined on a given interval $[x,y].$ Therefore, for the sake 
of simplicity, from this moment forth we shall work on the interval $[0,1].$

First we shall use Lemma \ref{mainl} to prove a result which extends the inequalities
from \cite{Szostok}.  

\begin{theorem}
\label{thlH}
 Let numbers $a_1,a_2,a_3,\alpha_1,\alpha_2,\alpha_3\in(0,1)$ satisfy
$a_1+a_2+a_3=1$ and $\alpha_1>\alpha_2>\alpha_3.$

Then the inequality 
\begin{equation}
\label{lH}
\sum_{i=1}^3a_if(\alpha_ix+(1-\alpha_i)y)\leq \frac{1}{y-x}\int_{x}^yf(t)
\end{equation}
is satisfied by all convex functions $f:[x,y]\to\eR$ if and only if we have 
\begin{equation}
\label{E1}
\sum_{i=1}^3a_i(1-\alpha_i)=\frac12
\end{equation}

 and one of the following conditions is satisfied

(i) $a_1\leq 1-\alpha_1$ and $a_1+a_2\geq 1-\alpha_3,$

(ii) $a_1\geq 1-\alpha_2$ and $a_1+a_2\geq 1-\alpha_3,$

(iii) $a_1\leq 1-\alpha_1$ and $a_1+a_2\leq 1-\alpha_2,$

(iv) $a_1\leq 1-\alpha_1, a_1+a_2\in(1-\alpha_2,1-\alpha_3)$ and $2\alpha_3\geq a_3,$

(v) $a_1\geq 1-\alpha_2,a_1+a_2<1-\alpha_3$ and $2\alpha_3\geq a_3,$

(vi) $a_1> 1-\alpha_1,a_1+a_2\leq 1-\alpha_2$ and $1-\alpha_1\geq\frac{a_1}{2},$

(vii) $a_1\in(1-\alpha_1,1-\alpha_2),$ $a_1+a_2\geq 1-\alpha_3,$ 
and $1-\alpha_1\geq\frac{a_1}{2}$

or

(viii) $a_1\in(1-\alpha_1,1-\alpha_2),$ $a_1+a_2\in(1-\alpha_2,1-\alpha_3),$ 
 $1-\alpha_1\geq\frac{a_1}{2}$ and $2a_1(1-\alpha_1)+2a_2(1-\alpha_2)\geq (a_1+a_2)^2.$
\end{theorem}
{\sc Proof} Let functions $F_1,F_2:\eR\to\eR$ be given by the following formulas
\begin{equation}
\label{F11}
F_1(t):=\left\{\begin{array}{ll}
0&t<1-\alpha_1\\
a_1&t\in[1-\alpha_1,1-\alpha_2)\\
a_1+a_2&t\in[1-\alpha_2,1-\alpha_3)\\
1&t\geq 1-\alpha_3
\end{array}\right.
\end{equation}
and
\begin{equation}
\label{F21}
F_2(t):=\left\{\begin{array}{ll}
0&t<0\\
t&t\in[0,1)\\
1&t\geq 1.
\end{array}\right.
\end{equation}
Observe that equality \rf{E1} gives us 
$$\int_{0}^1tdF_1(t)=\int_{0}^1tdF_2(t).$$
Further, it is easy to see that in cases $(i),(ii)$ and $(iii)$ the pair $(F_1,F_2)$ 
crosses exactly once and, consequently, inequality \rf{lH} follows from Ohlin lemma.

In case $(iv)$ pair $(F_1,F_2)$ crosses three times. Let $A_0,\dots,A_3$ be defined
so as in Lemma \ref{mainl}. In order to prove inequality \rf{lH} we must 
check that $A_0\geq A_1.$ However, since $A_0-A_1+A_2-A_3=0,$ we shall show that 
$A_2\leq A_3.$  We have 
$$A_2=\int_{a_1+a_2}^{1-\alpha_3}t-a_1-a_2dt=\frac{(1-\alpha_3-a_1-a_2)^2}{2}=\frac{a_3^2-2a_3\alpha_3+\alpha_3^2}{2}$$
and
$$A_3=\int_{1-\alpha_3}^1 1-tdt=\frac{\alpha_3^2}{2}.$$
This means that $A_2\leq A_3$ is equivalent to $2\alpha_3\geq a_3,$ as claimed.

We omit similar proofs in cases $(v),(vi)$ and $(vii)$ and we pass to the case $(vii).$
In this case the pair $(F_1,F_2)$ crosses five times. 
We have
$$A_0=\int_0^{1-\alpha_1}tdt=\frac{(1-\alpha_1)^2}{2}$$
and
$$A_1=\int_{1-\alpha_1}^{a_1}a_1-tdt=a_1(a_1-(1-\alpha_1))-\frac{a_1^2-(1-\alpha_1)^2}{2}=
\frac{[a_1-(1-\alpha_1)]^2}{2}.$$
This means that inequality $A_0\geq A_1$ is satisfied if and only if $1-\alpha_1\geq\frac{a_1}{2}.$

Further 
$$A_2=\int_{a_1}^{1-\alpha_2}t-a_1dt=\frac{(1-\alpha_2)^2-a_1^2}{2}-a_1(1-\alpha_2-a_1)$$
and 
$$A_3=\int_{1-\alpha_2}^{a_1+a_2}a_1+a_2-tdt=(a_1+a_2)(a_1+a_2-(1-\alpha_2))-\frac{(a_1+a_2)^2-(1-\alpha_2)^2}{2}$$
 therefore inequality $A_0+A_2\geq A_3+A_1$ is satisfied if and only if
$$(1-\alpha_1)^2+(1-\alpha_2-a_1)^2\geq (a_1-1-\alpha_1)^2+(a_1+a_2-1+\alpha_2)^2$$
which after some calculations gives us the last inequality from $(vii).$

Using  assertions (i) and (vii) of Theorem \ref{thlH}, it is easy to get the following example. 

\begin{example}
\label{pr1}
Let $x,y\in\eR,$ let  $\alpha\in(\frac12,1)$  and let $a,b\in(0,1)$ satisfy $2a+b=1.$ 
Then inequality 
\begin{equation}
\label{symlH}
af(\alpha x+(1-\alpha)y)+bf\left(\frac{x+y}{2}\right)+af((1-\alpha)x+\alpha y)\leq\frac{1}{y-x}\int_x^yf(t)dt
\end{equation}
is satisfied by all convex functions $f:[x,y]\to\eR$ if and only if  $a\leq 2-2\alpha.$
\end{example}
In the next theorem we shall obtain inequalities which extend the second of Hermite-Hadamard
inequalities.
\begin{theorem}
\label{thrH}
 Let numbers $a_1,a_2,a_3,a_4\in(0,1),\alpha_1,\alpha_2,\alpha_3,\alpha_4\in[0,1]$ satisfy
$a_1+a_2+a_3+a_4=1$ and $1=\alpha_1>\alpha_2>\alpha_3>\alpha_4=0.$

Then the inequality 
\begin{equation}
\label{rH}
\sum_{i=1}^4a_if(\alpha_ix+(1-\alpha_i)y)\geq \frac{1}{y-x}\int_{x}^yf(t)
\end{equation}
is satisfied by all convex functions $f:[x,y]\to\eR$ if and only if we have 
\begin{equation}
\label{E4}
\sum_{i=1}^4a_i(1-\alpha_i)=\frac12
\end{equation}
and one of the following conditions is satisfied:

(i) $a_1\geq1-\alpha_2$ and $a_1+a_2\geq1-\alpha_3,$

(ii) $a_1+a_2\leq 1-\alpha_2$ and $a_1+a_2+a_3\leq 1-\alpha_3,$ 

(iii) $1-\alpha_2\leq a_1$ and $1-\alpha_3\geq a_1+a_2+a_3,$

(iv) $1-\alpha_2\leq a_1,$ $1-\alpha_3\in(a_1+a_2,a_1+a_2+a_3)$ and $\alpha_3\leq 2a_4,$ 

(v) $1-\alpha_2\geq a_1+a_2,a_1+a_2+a_3>1-\alpha_3$ and $\alpha_3\leq 2a_4,$

(vi) $a_1<1-\alpha_2,$ $a_1+a_2\geq 1-\alpha_3$ and $2a_1+\alpha_2\geq 1,$

(vii) $a_1<1-\alpha_2,a_1+a_2>1-\alpha_2,a_1+a_2+a_3\leq 1-\alpha_3$ and $2a_1+\alpha_2\geq 1,$

(viii) 
$1-\alpha_2\in(a_1,a_1+a_2),1-\alpha_3\in(a_1+a_2,a_1+a_2+a_3),2a_1+\alpha_2
\geq 1$ and $2a_1(1-\alpha_3)+2a_2(\alpha_2-\alpha_3)\geq(1-\alpha_3)^2.$
\end{theorem}
{\sc Proof} 
Let function $F_1:\eR\to\eR$ be given by the following formula
\begin{equation}
\label{F12}
F_1(t):=\left\{\begin{array}{ll}
0&t<0\\
a_1&t\in[0,1-\alpha_1)\\
a_1+a_2&t\in[1-\alpha_1,1-\alpha_2)\\
a_1+a_2+a_3&t\in[1-\alpha_2,1)\\
1&t\geq 1
\end{array}\right.
\end{equation}
and let $F_2$ be given by \rf{F21}.
In view of \rf{E4}, we have 
$$\int_{0}^1F_1(t)dt=\int_{0}^1F_2(t)dt.$$
In cases $(i)-(iii)$ there is only one crossing point of $(F_2,F_1)$  
and our assertion is a consequence of Ohlin lemma.

In cases $(iv)-(vii)$ pair $(F_2,F_1)$ crosses three times and, therefore we have to use Lemma \ref{mainl}.
For example in case $(iv)$ \rf{rH} is satisfied by all convex functions $f$ if and only if $A_0\geq A_1.$ Further we know that 
$$A_0-A_1+A_2-A_3=0$$
thus $A_0\geq A_1$ is equivalent to $A_3\geq A_2.$ We clearly have
\begin{equation}
\label{A2}
\begin{split}
A_2=&\int_{1-\alpha_3}^{1-a_4}F_1(t)-F_2(t)dt=(\alpha_3-a_4)(1-a_4)-\frac{(1-a_4)^2-(1-\alpha_3)^2}{2}\\
=&(\alpha_3-a_4)\left(1-a_4+\frac{2-(\alpha_3+a_4)}{2}\right)
\end{split}
\end{equation}
and 
\begin{equation}
\label{A3}
A_3=\int_{1-a_4}^1 t-(1-a_4)dt=\frac{1-(1-a_4)^2}{2}-(1-a_4)a_4
\end{equation}
 i.e. $A_3\geq A_2$ is equivalent to $\alpha_3\leq 2a_4.$

We omit similar reasonings in cases $(v),(vi)$ and $(vii)$ and we pass to the most interesting case $(viii).$
In this case $(F_2,F_1)$ has 5 crossing points and, therefore, we must check that inequalities 
$$A_0\geq A_1\;\;\textrm{and}\;\;A_0-A_1+A_2\geq A_3$$ are equivalent to
respective inequalities of the condition $(viii).$
To this end we write
$$A_0=\int_{0}^{a_1}a_1-tdt=\frac{a_1^2}{2},$$
$$A_1=\int_{a_1}^{1-\alpha_1}t-(a_1+a_2)dt=\frac{(a_1+a_2-1+\alpha_1)^2}{2}$$
which means that $A_0\geq A_1$ if and only if $2a_1+\alpha_2\geq 1.$ 
Further $A_2$ and $A_3$ are given by formulas \rf{A2} and \rf{A3}. Thus
$A_0-A_1+A_2\geq A_3$ is equivalent to 
$$a_1^2+(a_1+a_2-(1-\alpha_2))^2\geq(1-\alpha_2-a_1)^2+(1-\alpha_3-a_1-a_2)^2$$
which yields
$$2a_1(1-\alpha_3)+2a_2(\alpha_2-\alpha_3)\geq(1-\alpha_3)^2.$$
\begin{remark} 
In this remark we use notations from papers \cite{BP} and \cite{EF}.
Results obtained in \cite{Szostok} allowed to get the inequalities 
$$\frac{f(x)+f(a+b-x)}{2}\leq\frac{1}{b-a}\int_{a}^bf(t)dt,\;x\in\left[\frac{3a+b}{4},\frac{a+b}{2}
\right]$$
which was proved in \cite{BP}
and
$$f\left(\frac{x+y}{2}\right)\leq\lambda f\left(\frac{\lambda b+(2-\lambda)a}{2}\right)
+(1-\lambda)f\left(\frac{(1+\lambda) b+(1-\lambda)a}{2}\right)\leq$$
$$\frac{1}{b-a}\int_a^bf(t)dt\leq\frac12(f(\lambda b+(1-\lambda)a)+
\lambda f(a)+(1-\lambda) f(b))\leq$$$$\frac12 (f(x)+f(y)),\;\;\lambda\in[0,1]$$
from \cite{EF}. However, it was not possible to apply these results to prove the inequality
\begin{equation}
\label{BP1}
\frac{f(a)+f(x)+f(a+b-x)+f(b)}{4}\geq\frac{1}{y-x}\int_{a}^bf(t)dt,\;\;x\in\left[a,\frac{a+b}{2}\right]
\end{equation}
(see \cite{BP}).
Now, inequality \rf{BP1} is an easy consequence of Theorem \ref{thrH}.
\end{remark}
It is interesting to note that it is easy to get a more general inequality than \rf{BP1}.
Namely, using assertions $(ii)$ and $(vii)$ of Theorem \ref{thrH} 
we get the following example.
\begin{example}
Let $x,y\in\eR,$ let $a,b\in(0,1),\alpha\in(\frac12,1)$ and let $a,b$ satisfy $2a+2b=1.$ 
Then inequality 
$$af(x)+bf(\alpha x+(1-\alpha)y)+bf\left((1-\alpha)x+\alpha y\right)+af(y)\geq\frac{1}{y-x}\int_x^yf(t)dt$$
is satisfied by all convex functions $f:[x,y]\to\eR$ if and only if 
$a\geq\frac{1-\alpha}{2}.$
\end{example}

In the next part of this section we show that the same tools may be used to obtain some inequalities between quadrature operators which do not involve the integral mean.  
\begin{theorem}
\label{thqo}
 Let $a,\alpha_1,\alpha_2,\beta\in(0,1)$ and let $b_1,b_2,b_3\in(0,1)$ satisfy
$b_1+b_2+b_3=1.$

Then the inequality 
\begin{equation}
\label{inqo}
af(\alpha_1x+(1-\alpha_1)y)+(1-a)f(\alpha_2x+(1-\alpha_2)y)\leq b_1f(x)+b_2f(\beta x+(1-\beta)y)+b_3f(y)
\end{equation}
is satisfied by all convex functions $f:[x,y]\to\eR$ if and only if we have 
\begin{equation}
\label{E}
b_2(1-\beta)+b_3=a(1-\alpha_1)+(1-a)(1-\alpha_2)
\end{equation}

 and one of the following conditions is satisfied:

(i) $a\leq b_1,$

(ii) $a\geq b_1+b_2,$

(iii) $\alpha_2\geq\beta$

or

(iv) $a\in(b_1,b_1+b_2),$ $\alpha_2<\beta$ and $(1-\alpha_1)b_1\geq(\alpha_1-\beta)(a-b_1).$

\end{theorem}
{\sc Proof}
Define 
\begin{equation}
\label{F1}
F_1(t):=\left\{\begin{array}{ll}
0&t<1-\alpha_1\\
a&t\in[1-\alpha_1,1-\alpha_2)\\
1&t\geq 1-\alpha_2
\end{array}\right.
\end{equation}
and
\begin{equation}
\label{F2}
F_2(t):=\left\{\begin{array}{ll}
0&t<0\\
b_1&t\in[0,1-\beta)\\
b_1+b_2&t\in[1-\beta,1)\\
1&t\geq 1
\end{array}\right.
\end{equation}
Then inequality \rf{inqo} takes form
\begin{equation}
\label{inqoI}
\int_x^y fdF_1\leq\int_x^y fdF_2.
\end{equation}
Moreover \rf{E} means that
$$\int_x^y tdF_1(t)=\int_x^y tdF_2(t).$$
It is easy to see that in cases (i),(ii) and (iii) $F_1$ and $F_2$ cross exactly once 
and, therefore, our assertion follows from Ohlin lemma.

Now we pass to the most interesting case $(iv).$ In this case functions $F_1$ and $F_2$ 
cross three times, at points: $1-\alpha_1,1-\beta$ and $1-\alpha_2.$ Thus
$$A_0=\int_0^{1-\alpha_1}tdF_2(t)=b_1(1-\alpha_1)$$     
$$A_1=\int_{1-\alpha_1}^{1-\beta}tdF_2(t)=(\alpha_1-\beta)(a-b_1).$$
In view of Lemma \ref{mainl}, we know that inequality \rf{inqoI}
is satisfied (for all convex functions $f$) if and only if $A_0\geq A_1$
which ends the proof.
 
Now, using this theorem, we shall present positive and negative examples of inequalities of the type  \rf{inqo}.
\begin{example} Let $\alpha\in\left(\frac12,1\right)$
Inequality 
$$\frac{f(\alpha x+(1-\alpha)y)+f((1-\alpha)x+\alpha y)}{2}\leq\frac{f(x)+f\left(\frac{x+y}{2}\right)+f(y)}{3}$$
is satisfied by all convex functions $f:[x,y]\to\eR$ if and only if $\alpha\leq\frac56.$
\end{example}

\begin{example} Let $\alpha\in\left(\frac12,1\right)$
Inequality 
$$\frac{f(\alpha x+(1-\alpha)y)+f((1-\alpha)x+\alpha y)}{2}\leq\frac16f(x)+\frac{2}{3}f\left(\frac{x+y}{2}\right)+\frac16f(y)$$
is satisfied by all convex functions $f:[x,y]\to\eR$ if and only if $\alpha\leq\frac23.$
\end{example}

\begin{remark}
It is clear that our results  may be 
extended to cover longer expressions but the calculations would become 
very complicated  and the best way to deal with it may be to write 
a computer program  which would check if in concrete cases inequalities of such types are satisfied.
\end{remark}

\begin{remark} As it is known from the paper \cite{BesPal}, if a continuous function
satisfies inequalities  
of the type which we have considered then such function must be convex.

Therefore inequalities obtained in this paper characterize convex functions 
(in the class of continuous functions).
\end{remark}


\begin{thebibliography}{20}
\bibitem{BesPal}
 M. Bessenyei, Zs. P\'{a}les, 
{\it  Characterization of higher order monotonicity via
integral inequalities,} Proc. Roy. Soc. Edinburgh Sect. A 140 (2010), no. 4, 723--736. 

\bibitem{DP}
        {S.S. Dragomir, C.E.M. Pearce}, 
        {\it Selected Topics on Hermite--Hadamard Inequalities and Applications}, 
         2000.

\bibitem{EF}A. El Farissi, {\it Simple proof and refinement of Hermite-Hadamard inequality} J. Math. Inequal.  4  (2010),  no. 3, 365–369.

 \bibitem{FPPV} I.Franjic, J. Pecaric, I. Peric, A. Vukelic,{\it Euler integral identity, quadrature formulae and error estimations (from the point of view of inequality theory)} Monographs in Inequalities, 2. ELEMENT, Zagreb, 2011.

\bibitem{BP} 
M. Klaricic Bakula, J. Pecaric {\it Generalized Hadamard's inequalities based on general Euler 4-point formulae} ANZIAM J.  48  (2007),  no. 3, 387–404.

\bibitem{LS}V.I. Levin, S.B. Ste\v{c}kin; 
{\it Inequalities.} 
Amer. Math. Soc. Transl. (2)  14  1960 1--29.

\bibitem{NP}C. P. Niculescu, L-E. Persson, 
Convex functions and their applications. 
A contemporary approach. CMS Books in Mathematics/Ouvrages de Math\'ematiques de la SMC, 23. Springer, New York, 2006.

\bibitem{Ohlin}
        {J. Ohlin},
        {\it On a class of measures of dispersion with application to optimal reinsurance},
        ASTIN Bulletin {\bf 5} (1969), 249--66.

\bibitem{OS1} A. Olbry\'s, T. Szostok, {\it Inequalities of the Hermite-Hadamard type involving numerical differentiation formulas,} manuscript

\bibitem{Rajba}{T. Rajba},
         {\it  On The Ohlin lemma for Hermite-Hadamard-Fejer type inequalities}, Math. Ineq. Appl. {\bf 17}, Number 2 (2014), 557--571.

\bibitem{Szostok}{T. Szostok} {\it Ohlin's lemma and some inequalities of the Hermite-Hadamard type,} Aequationes
Math., DOI:10.1007/s00010-014-0286-2
\end{thebibliography}
\end{document}